\documentclass[12pt]{article} 
\usepackage{verbatim,amsmath,amssymb}

\newcommand{\rf}[1]{(\ref{#1})}

\newcommand{\numer}[1]{\label{#1}}
\newcommand{\bibl}[1]{\bibitem{#1}}

\newcommand{\cd}{{\cal D}}
\newcommand{\cc}{{\cal C}}

\newcommand{\cu}{{\cal U}}
\newcommand{\ch}{{\cal H}}

\newcommand{\ct}{{\cal T}}

\usepackage{epsfig}
\usepackage{graphicx}
\usepackage{amscd}
\input math111.def
\input buch.def
\numberwithin{equation}{section}
\newif\ifmarglab
\makeatletter
\def
\label#1{\@bsphack\ifmarglab\marginpar{LAB:#1}\fi\if@filesw {\let\thepage\relax
   \def\protect{\noexpand\noexpand\noexpand}%
   \edef\@tempa{\write\@auxout{\string
      \newlabel{#1}{{\@currentlabel}{\thepage}}}}%
   \expandafter}\@tempa
   \if@nobreak \ifvmode\nobreak\fi\fi\fi\@esphack}
\makeatother
\begin{document}
\title{Combinatorics of Dyadic Intervals: Consistent Colourings 
}
\author{Anna Kamont, Paul F. X . M\"uller}
%\date{November 29, 2013}
\maketitle

\begin{abstract}
 
We study the problem of consistent and homogeneous colourings for increasing families of dyadic intervals. 
We determine when this problem can be solved and when not.

\end{abstract}

MSC 2010: 05A18,  91A05

Key words: Coloured dyadic intervals, consistent and homogeneous colouring, 2-person games

\section{\numer{c.s1} Introduction}

Combinatorics of coloured dyadic intervals refers to a set of techniques created for the study of operators
defined through their action on the Haar system. We refer to 
the treatment of averaging projections by P.W.~Jones \cite{jones-85},
the proof of the vector-valued $T(1)$ theorem by T~ Figiel \cite{figiel-1, figiel-2}, the use of the stripe operators
in J.~Lee, P.F.X.~M\"uller, S.~M\"uller \cite{lmm}, and the study of rearrangement operators on $L^p$ spaces,
 P.F.X.~M\"uller \cite{pfxm-bmo}, K.~Smela \cite{smela},
A.~Kamont, P.F.X.~M\"uller \cite{ak-pm}.

Here we study a very natural colouring problem on dyadic trees. We start out with a coloured collection of dyadic intervals $\cc$,
where we assume that the colours are distributed homogeneously over $\cc$.
Given any collection $\ch$ containing $\cc$ 
we ask if there exists an equally homogeneous colouring of $\ch$ that preserves the colours of $\cc$
(consistent colouring of $\ch$). The nature of this problem depends
very much on what we agree to call a homogeneous distribution of colours.
Our choice of homogeneity is very restrictive, and consequently in working on the problem of
consistent colouring we encountered delicate combinatorial questions.

\bigskip

Let $\cd$ denote the collection of dyadic intervals in the unit interval $[0,1]$, and let
$$
\cd_j = \{ I \in \cd: |I| = 2^{-j}\}.
$$
We consider a large collection $\cc \subset \cd_j$.
We assume that the intervals in $\cc$ are painted with $d$ distinct colours, giving rise to a decomposition
$$
\cc = \cc_1 \cup \ldots \cup \cc_d.
$$
It is intuitively clear what it means that the colours $\{1, \ldots, d\}$ are homogeneously distributed among the intervals of $\cc$.
For instance, we would demand that there exists $\eta > 0$ so that
\begin{equation}
\numer{z.1}
\eta \max_{1 \leq i \leq d} |\cc_i| \leq \min_{1 \leq i \leq d} |\cc_i|,
\end{equation}
where $|\cc_i|$ denotes the cardinality of the collection $\cc_i$. A much stronger measure of homogeneity
arises when we ask for \rf{z.1} to hold over the prespecified collection of testing intervals
$$
\ct = \{J \in \cd: |J| > 2^{-j}\}.
$$
Specifically, if $|\cc \cap L| > d$, we would demand that there exists $\eta >0$ so that
\begin{equation}
\numer{z.2}
\eta \max_{1 \leq i \leq d} |\cc_i \cap L | \leq \min_{1 \leq i \leq d} |\cc_i \cap L | 
\quad \hbox{ for each }\quad  L \in \ct,
\end{equation}
where
$$
\cc_i \cap L =\{I \in \cc_i: I \subset L\}.
$$

\smallskip

We use an additional rule to express homogeneity with respect to testing intervals that satisfy $|\cc \cap L|\leq d$.
The necessity of such a rule arises 
from the fact that the cardinalities $|\cc_i \cap L|$  take values in $\bN \cup \{0\}$, hence if
$|\cc \cap L| < d$, then \rf{z.2} has to fail.
 Thus, there are two regimes -- high cardinality and low cardinality of $\cc \cap L$, and the transistion
arises at $|\cc \cap L| = d$. The following definition contains the homogeneity conditions for both regimes, and 
 it addresses
the discrete nature of our gauge functions
$$\cc_i \to |\cc_i \cap L|, \quad L \in \ct.$$

\begin{defi}
\numer{c.d1}
Let $\cc \subset \cd_j $, and fix $d \in \bN$, $0<\eta \leq {1 \over 2}$. 
Let   $\cc = \cc_1 \cup \ldots \cup \cc_d$ be some decomposition of $\cc$.
This decomposition is called $(\eta,d)$-homogeneous  colouring of $\cC$ if for each $L \in \cd$, $|L| \geq {1 \over 2^j}$
one of the following holds: 
\begin{itemize}
\item[]
{\bf \em Either } $|\cc \cap L| >d$, and then
\begin{equation}
\numer{hom2}
\eta \max_{1 \leq i \leq d} |\cc_i \cap L| \leq  \min_{1 \leq i \leq d} |\cc_i \cap L|,
\end{equation}
\item[] 
{\bf \em Or else} $|\cc \cap L| \leq d$, and then 
\begin{equation}
\numer{hom1}
|\cc_i \cap L|  \leq 1 \quad \text{ for each } \quad 1 \leq i \leq d,
\end{equation}

\end{itemize}
\end{defi}

\paragraph{Remark.} We remark that 
for each (uncoloured) $\cc \subset \cd_j$, $d\in \bN$ and $\eta = {1 \over 2}$ there is always 
a $(\eta,d)$-homogeneous colouring that can be
obtained as follows: Enumerate the intervals in $\cc$ from left to right, 
and simply put %for $1 \leq r \leq d$
\begin{equation}
\numer{mod.d}
\cc_r = \{\Gamma_l \in \cc: l = r \hbox{ mod } d\},\quad\quad 1\le r \leq d.
\end{equation}
Later, we refer to such a colouring as a {\em colouring modulo $d$}.

\bigskip

The use of this colouring rule -- applied to intervals of equal length -- appeard in a context similar 
to ours in \cite[p. 200]{jones-80}, see also \cite[p. 359]{garnett} and \cite[p. 199]{pfxm}.

\paragraph{The problem of consistent colouring.} 
The problem we treat in this paper is 
the following. We are given two disjoint collections $\cc, \cu \subset \cd_j$.
Assume that the collection $\cc$ is coloured, that is, it is given
an $(\eta,d)$-homogeneous  colouring
$$
\cc = \cc_1 \cup \ldots \cup \cc_d.
$$
The collection $\cu$ consists of uncoloured intervals.
We would like to colour the intervals in $\cu$ with the same colours $\{1 , \ldots , d\}$, that is to decompose
$\cu$ as 
$$
\cu = \cu_1 \cup \ldots \cup \cu_d
$$
in such a way that the union $\ch = \cc \cup \cu$ has an $(\eta,d)$-homogeneous    colouring given by
$$
\ch = \ch_1 \cup \ldots \cup \ch_d, \quad
\hbox{ where } \quad
\ch_i = \cc_i \cup \cu_i \quad  \hbox{ for } \quad  1 \leq i \leq d.
$$
That is, we want to obtain an $(\eta,d)$-homogeneous colouring of $\ch \supset \cc$ keeping the pre-existing 
$(\eta,d)$-homogeneous  colouring of $\cc$.

We refer to this question as to {\em the problem of finding a colouring of $\ch$  consistent with existing colouring of $\cc$.}
Our treatment of this problem is as follows:

\begin{enumerate}
\item We isolate a condition on $\cU$ and $\cC$ (previsibility; see Definition~\ref{c.d2})
implying
that the   problem of consistent colouring for $\ch = \cc \cup \cU  $ has a solution. 
See Theorem~\ref{c.t1}.  
\item We give  examples where the  
problem of consistent colouring for $\ch = \cc \cup \cU  $
has just one solution. 
Moreover, we give examples (of $\cC$, its decomposition $\{\cC_i\}$ and $\cU$) for which
 the 
problem of consistent colouring for $\ch = \cc \cup \cU  $ 
  does not have a solution. See Proposition~\ref{e.s1}.  
\item In Section \ref{game} we reformulate the problem of consistent colouring as a two-person game. Our results -- Theorem \ref{c.t1}
and Proposition \ref{e.s1} -- translate into winning strategies for the respective players.
\end{enumerate}

For the appearence of succesive colourings
 of dyadic intervals in the context of averaging projections see \cite[p. 871-875]{jones-85}.
In \cite{ak-pm} we constructed supporting trees for some rearrangement operators 
and thereby proved their boundedness on vector valued $L_p$ spaces.
Initially, our approach to defining the supporting trees was by inclusion-exclusion principles and consistent colourings as studied 
in the present paper.

\section{\numer{c.s2} Constructing a consistent colouring}

In the following   we isolate a criterion which guarantees the existence of consistent colouring.
To formulate this criterion, we use a dyadic interval $L \in \cd$ together with its immediate dyadic successors $L',L''$, i.e.
intervals $L',L'' \in \cd$ such that $L = L' \cup L''$ and $|L'| = |L''| = {1 \over  2 } |L|$.

The problem of consistent colouring lead us to the following  condition:

\smallskip 

\noindent
We are given disjoint collections $\cc, \cu \subset \cd_j$
and $d \in \bN$. We say
 that the pair 
 $(\cc,\cu)$ is {\em $d$-previsible} if with $\ch = \cc \cup \cu$, the conditions
 $$
 |\ch \cap  L''| \geq d, \quad \cc \cap L'' \neq \emptyset,  \quad \cu \cap L'' \neq \emptyset
 $$
 imply 
 $$
 |\ch \cap L'| \geq d.
 $$
 To facilitate precise reference in the course of our argument below, we encode this notion of $d$-previsibility 
 in the following -- equivalent -- definition.

\begin{defi}
\numer{c.d2}
Let $\cc, \cu \subset \cd_j$, $\cc \cap \cu = \emptyset$. Let $d \in \bN$.
The pair of collections $(\cc, \cu)$ is called $d$-previsible 
 if for every $L \in \cd$ with $|L| \geq {1 \over 2^{j-1}}$ and its dyadic succesors $L', L''$, 
 the following holds: 
 $$
 |(\cu\cup \cc) \cap L'| < d \quad \hbox{and} \quad  |(\cu\cup \cc) \cap L''| \geq d \quad 
 \hbox{implies} \quad  \cu \cap  L''=\emptyset \quad \hbox{or} \quad \cc \cap  L''=\emptyset .
 $$
\end{defi}
Now, we have the following Theorem \ref{c.t1} which gives a sufficient condition for existence of consistent colourings.

\begin{theor}
\numer{c.t1}
Fix $d \in \bN$ and $\eta$, $0<\eta \leq {1 \over 2}$.
Let $\cc \subset \cd_j$, and let $\{\cc_i, 1 \leq i \leq d\}$ be a fixed $(\eta,d)$-homogeneous  colouring of $\cc$.
Let $\cu \subset \cd_j$ be such that the pair $(\cc, \cu)$ is $d$-previsible.
%with respect to $\cc$. 
Then there is a colouring
$\{\cu_i, 1 \leq i \leq d\}$ of  
$\cu$ such that $\{\ch_i = \cc_i \cup \cu_i, 1 \leq i \leq d\}$ is an  
$(\eta,d)$-homogeneous  colouring of $\ch = \cc\cup\cu$.
\end{theor}

\paragraph{Remark.} Ones first attempt to prove Theorem \ref{c.t1} by an inductive argument would be the following:
Find first $(\eta,d)$-homogeneous colouring of $\ch \cap K$ for $K \in \cd_{j - \alpha}$, where $2^\alpha \leq d < 2^{\alpha +1}$.
Then carry over these colourings -- inductively and backwards in time -- to 
larger collections $\ch \cap L$, $L \in \cd_s$ with $s > j-\alpha$ as follows:
Assume that for a dyadic interval $L$ with successors $L', L''$, the separate $(\eta,d)$-homogeneous colourings of
$\ch \cap L'$ and  $\ch \cap L''$ are fixed. Then check that the union of these colourings gives an $(\eta,d)$-homogeneous colouring
of the union $\ch \cap L = (\ch \cap L') \cup (\ch \cap L'')$. 
If this procedure would work, at stage $s$,
we would have produced an $(\eta,d)$-homogeneous colouring of $\ch \cap K$ for each $K \in \cd$ with $|K| \leq 2^{-s}$.
However, such a deterministic approach cannot work, as the following example shows.

\paragraph{Example.} Take $d = 3 r$, $r \in \bN$, and fix a dyadic interval $L$ with successors $L', L''$.
Assume that collections $\cc, \cu$ are such that $| \cc \cap L'| = |\cc \cap  L''| = r$ and 
$|\cu\cap L'| = |\cu \cap  L''| = r$. Take an $(\eta,d)$-homogeneous colouring of $\cc \cap L$ such that
$$
|\cc_i \cap L'| =1 \quad \hbox{ for } \quad 1 \leq i \leq r 
\quad\hbox{ and } \quad |\cc_i \cap L'|= 0  \quad\hbox{ for }  \quad r+1 \leq i \leq 3r,
$$
$$
 |\cc_i \cap L''|=1 \quad\hbox{ for } \quad r+1 \leq i \leq 2r, \quad |\cc_i \cap L''| = 0 \quad \hbox{ for }\quad 1 \leq i \leq r
\quad \hbox{ and } \quad 2r+1 \leq i \leq 3r.
$$
Note that then we have 
$$
 |\cc_i \cap L|=1 \quad \hbox{ for }\quad 1 \leq i \leq 2r \quad \hbox{ and }\quad |\cc_i \cap L| = 0 
 \quad\hbox{ for } \quad 2r+1 \leq i \leq 3r.
$$
Next, we can choose colourings of
$\cu \cap L'$ and $\cu \cap L''$ such that 
$$
|\cu_i \cap L'| =1 \quad \hbox{ for }\quad  r+1 \leq i \leq 2r, \quad  |\cu_i \cap L'|= 0 \quad \hbox{ for } \quad 1 \leq i \leq r
\quad \hbox{ and } \quad 2r+1 \leq i \leq 3r,
$$
$$
 |\cu_i \cap L''|=1 \quad \hbox{ for } \quad 1 \leq i \leq r \quad \hbox{ and } \quad  |\cu_i \cap L''|= 0 \quad
 \hbox{ for } \quad r+1 \leq i \leq 3r.
$$
Then, knowing the colouring of $\cc \cap K$ and $\cu \cap K$ for $K = L', L''$ with $\ch = \cc \cup \cu$ we have
$$
| \ch_i \cap  K| =1\quad \hbox{ for }\quad 1 \leq i \leq 2r \quad \hbox{ and } \quad  | \ch_i \cap  K|= 0 
\quad  \hbox{ for } \quad 2r+1 \leq i \leq 3r,
$$
so we have separate $(\eta,d)$-homogeneous colourings of $\ch \cap L'$ and $\ch \cap L''$.
However, by taking the union of these colourings we get a colouring of $\ch \cap L$ such that 
$$
| \ch_i \cap  L| =2  \quad\hbox{ for } \quad 1 \leq i \leq 2r \quad\hbox{ and } \quad
| \ch_i \cap  L| = 0 \quad\hbox{ for } \quad 2r+1 \leq i \leq 3r,
$$
which is not $(\eta,d)$-homogeneous.

Note however that the above example does not contradict the assertion of Theorem \ref{c.t1}. 
In fact -- given $\cc$ and its colouring as above -- it is quite easy to obtain colouring 
of $\cu$ so that the conclusion of Theorem \ref{c.t1} holds.
For  $K = L', L''$, we put  
$$
| \cu_i \cap  K| =0 \quad\hbox{ for } \quad 1 \leq i \leq 2r \quad \hbox{ and } 
\quad | \cu_i \cap  K| = 1 \quad \hbox{ for } \quad 2r+1 \leq i \leq 3r.
$$
Taking the union with the colouring of $\cc \cap K$, we find that $0 \leq | \ch_i \cap  K| \leq 1$, so we have an 
$(\eta,d)$-homogeneous colourings of $\ch \cap K$, $K = L', L''$. Finally,
$$
 | \ch_i \cap  L|=1 \quad\hbox{ for } \quad 1 \leq i \leq 2r \quad \hbox{ and } \quad | \ch_i \cap  L|= 2 \quad\hbox{ for }
 \quad 2r+1 \leq i \leq 3r,
$$
so we have an $({1 \over 2}, d)$-homogeneous colouring of $\ch \cap L$.\endproof

\bigskip

In response to  these examples, we introduced a stopping time argument -- running backwards in time -- that
produces the $(\eta,d)$-homogeneous colouring of Theorem \ref{c.t1}. 
At stage $s$ of our inductive argument, we will produce consistent $(\eta,d)$-homogeneous colourings of collections $\ch \cap K$
for $K \in \cd $ with $K \leq 2^{-s}$ provided that $K$ satisfies 
$$
| \ch \cap K| \geq d \hbox{ and } \cc \cap K \neq \emptyset.
$$

\paragraph{Proof of Theorem \ref{c.t1}.}
 We are going to define colouring of $\cu$ by an inductive argument.
Let $\alpha$ be such that $2^{\alpha }  \leq d < 2^{\alpha+1}$. 
Let us observe that if 
${1 \over 2^j} \leq |L| \leq {1 \over 2^{j-\alpha}}$, then $|\cH\cap L| \leq 2^\alpha \leq d$.
Thus, if the homogeneity
conditions \eqref{hom1} respectively \eqref{hom2}  are satisfied for $L \in \cd$
with $|L| \geq {1 \over 2^{j-\alpha}}$, then they are satisfied for each $L \in \cd$ with $|L| \geq {1 \over 2^j}$.
Therefore, in our procedure of colouring $\cu$  we consider only $L \in \cd_k$ with $k \leq j-\alpha$.

\medskip

The inductive argument is used to prove the following statement at each stage $s$, $j - \alpha \geq s \geq 0$: 

\medskip

\noindent
{\bf Inductive hypothesis at the stage $s$:} Let $K \in \cd_{s}.$ 
 If $|\ch\cap K| < d$ or $\cc \cap K = \emptyset$, then intervals in $\cu \cap K$  are still uncoloured. 
If $|\ch \cap K| \geq d$ and $\cc \cap K \neq \emptyset$, then all intervals in $\cu \cap K$ 
are coloured, 
and the  colouring of $\ch \cap K$ is $(\eta,d)$-homogeneous; 
as $|\ch \cap K| \geq d$, this means that
 $|\ch_i \cap K| \geq 1$ and
\begin{equation}
\numer{(b3)}
\eta \max_{1 \leq i \leq d} |\ch_i \cap K| \leq \min_{1 \leq i \leq d} |\ch_i \cap K|.
\end{equation}

\medskip

\noindent
{\bf I. The start of the induction.} Let $L \in \cd_{j-\alpha}$. Then either  $|\ch\cap L|<d$ or  $|\ch \cap L|=d$.

\smallskip

\noindent
{\bf I.1.} If $|\ch \cap L|<d$ or $\cc \cap  L = \emptyset$, then intervals in $\cu \cap L$  are left uncoloured.

\smallskip

\noindent
{\bf I.2.} If $|\ch \cap L|=d$ and $\cc \cap  L \neq \emptyset$, then also $| \cc \cap L| \leq d$, which implies that
$|\cc_i \cap L| \leq 1$ for each $1 \leq i \leq d$. In such case it is possible to colour
intervals in $\cu \cap L$  so that $|\ch_i \cap  L| = 1$
for each $i$, $1 \leq i \leq d$.

\medskip

\noindent
{\bf II. The inductive step.}   Let $\nu < j - \alpha$. 
The inductive assumption states that  there is a colouring at stage $\nu+1$. We need to prove that there is a colouring
at stage $\nu$. For this take
$L \in \cd_\nu$.  Then $L = L' \cup L''$
with $L' , L'' \in \cd_{\nu+1}$, $\nu+1 \leq j - \alpha$.
Each interval in $\cc \cap L$ or $\cu \cap L$  is included in $L'$ or
$L''$, so we have
\begin{eqnarray*}
 | \cc \cap  L| &  = & | \cc \cap  L'| + | \cc \cap  L''|,
\\
| \cu \cap  L| &  = & | \cu \cap  L'| + | \cu \cap  L''|,
\\
| \ch \cap  L| &  = & | \ch \cap  L'| + | \ch \cap  L''|.
\end{eqnarray*}

Then we have two main cases:

\medskip

\noindent
{\bf II.1.} $| \ch \cap L| < d$ or $ \cc \cap L  =\emptyset$. 
If  $| \ch \cap L| < d$  then also 
$| \ch \cap L'|, | \ch \cap L''| < d$ .
 If  $ \cc \cap L  =\emptyset$, then also $ \cc \cap L'  =\emptyset$ and $ \cc \cap L''  =\emptyset$.
In both cases, by the induction hypothesis,
 intervals in both  $\cu \cap L'$ and $\cu \cap L''$
are  uncoloured. 

If $\nu >0$, then leave intervals in $\cu \cap L$
 still uncoloured.

If $\nu=0$, then $L=[0,1]$, and the induction ends. This means that 
$|\ch| < d$ or $\cc = \emptyset$. 
If $|\ch|<d$, then  it is enough to assign elements of $\cu$ to
colours different from colours of elements of  $\cc$.
If $\cc = \emptyset$,  then it is enough to colour $\ch$ by modulo $d$ method, see  \rf{mod.d}.

\medskip

\noindent
{\bf II.2.} $| \ch \cap  L| \geq d$ and $\cc \cap L \neq \emptyset$. 
It follows that at least one of collections $\cc \cap L', \cc \cap L''$ must be nonempty.
Now we separate next two subcases: 

\smallskip

\noindent{\bf II.2.A.}  $| \ch \cap  L| \geq d$, $\cc \cap L \neq \emptyset$
and both $\cc \cap L' \neq \emptyset$, $\cc \cap L'' \neq \emptyset$.

\smallskip

\noindent{\bf II.2.B.} $| \ch \cap  L| \geq d$, $\cc \cap L \neq \emptyset$,
$\cc \cap L' = \emptyset$, but $\cc \cap L'' \neq \emptyset$.
(The case $\cc \cap L' \neq \emptyset$ and $\cc \cap L'' = \emptyset$
is symmetric to this one, 
and there is no need to treat it separately.)

\smallskip

We first treat the case II.2.A. Then we treat the case II.2.B.
Both these cases have their own subcases.  
The case II.2.B is treated by reducing its subcases to appropriate subcases of II.2.A.

\smallskip

\noindent
{\bf II.2.A.1.} $| \ch \cap  L'| \geq d$ and $|\ch \cap  L''| \geq
d$. 
Recall that $\cc \cap L' \neq \emptyset$ and $\cc \cap L'' \neq \emptyset$.
Then by induction hypothesis all intervals in $\cu \cap L' $  and in $\cu \cap L'' $ 
are already coloured, i.e. all intervals in $\cu \cap L$
  are coloured. Moreover, by \rf{(b3)}, for each $1 \leq i,k \leq d $
$$
\eta | \ch_i \cap L| = \eta  | \ch_i \cap L'| + \eta |\ch_i \cap L''|
\leq | \ch_k \cap L'| +  |\ch_k \cap L''|  =  | \ch_k \cap L| .
$$
Of course, we have also $ | \ch_i \cap L| \geq 1$.

\medskip

\noindent
{\bf II.2.A.2.} $|\ch \cap L'| < d$ and $|\ch \cap  L''| <
d$. Then by induction hypothesis all intervals in $\cu \cap L'$ 
and in  $\cu \cap L''$ are uncoloured, but the intervals in $\cc \cap L$ carry
 their colours. 

Now, we need to colour all  intervals in $\cu \cap L = (\cu \cap L') \cup (\cu \cap L'')$. 
To simplify notation, let
$$
m=|\cc \cap  L'|, \quad n = |\cc \cap  L''|, \quad x=|\cu \cap L'|,
\quad
y=|\cu \cap L''|.
$$
 We have 
$$
0 \leq m, n \leq d-1,  \quad 0 \leq m+x,n+y
\leq d-1\quad \hbox{ and }\quad d \leq |\ch \cap  L| = m+x+n+y \leq 2(d-1).
$$

First consider the case $m+n <d$. Then
 $$0 \leq |\cc_i \cap L| \leq
1 \quad \hbox{ for each } \quad i, \quad 1 \leq i \leq d.
$$ 
For simplicity, assume that intervals in $\cc \cap L'$ 
have colours $1, \ldots, m$, and intervals  in $\cc \cap L''$ have colours $m+1, \ldots , m+n$.
Now, we colour intervals in $\cu \cap L$.
First, colour intervals in $\cu \cap L'$ using 
colours $m+n+1, \ldots, d$, and then, if necessary (i.e. $x > d - (m+n)$), continuing with $x
- (d - (m+n))$
colours from  $m+1, \ldots, m+n$; since $m+x <d$, in this way we assign
colours to all intervals in  $\cu\cap L'$.
 Next, we assign colours to
intervals in  $\cu\cap L''$. If $m+n+x<d$, then assign
first colours $m+n+x+1, \ldots, d$, then continue with colours $1,
\ldots, m$, and then if necessary with colours $m+n+1, \ldots, m+n+x$. 
If $m+n+x \geq d$, then just choose $y$ different colours from $1,
\ldots,m$ and $m+n+1, \ldots, d$. With such colouring of intervals
 in $\cu \cap L'$  and in $\cu \cap L''$ we find that 
both
$$|\ch_i \cap L'| \leq 1 \quad \hbox{ and } \quad
  |\ch_i \cap L''| \leq 1, \quad 1 \leq i \leq d.
  $$
 This implies that for each $K \subset L'$ or $K\subset
L''$ we have $| \ch \cap K| < d$ and $|\ch_i\cap  K| \leq 1$.
Moreover, we get $1 \leq | \ch_i \cap L| \leq 2$, which implies
$$
\eta \max_{1 \leq i \leq d} |\ch_i \cap L| \leq {1 \over 2} \max_{1 \leq i \leq d}  |\ch_i \cap L|
\leq \min_{1 \leq i \leq d}  |\ch_i \cap L|.
$$

It remains to consider the case $m+n \geq d$. 
Then the homogeneity assumption on the decomposition of $\cc$ -- \eqref{hom1} for $L',L''$ and \eqref{hom2} for $L$ --  implies
$$
1 \leq |\cc_i \cap L|  = |\cc_i \cap L'| + |\cc_i \cap L''| \leq 2, \quad 1 \leq i \leq d.
$$
 For simplicity, assume
that intervals in $\cc \cap L'$  have colours $1, \ldots,
m$ and intervals in $\cc \cap L''$  have colours $m+1, \ldots, d$ and
$1, \ldots, m+n-d$ (note that $m+n-d <m$, since by assumption $n <d$). To colour intervals in $\cu \cap L'$ 
choose $x$ colours from $m+1, \ldots, d$. To colour intervals in
$\cu \cap L''$  choose $y$ colours from $m+n-d+1, \ldots , m$.
This is possible since $m+x < d$ and $n+y<d$. Observe that in this way
we get
$$
0 \leq |\ch_i \cap L'|, |\ch_i \cap L''| \leq 1 \quad \hbox{ and } \quad 1
\leq |\ch_i \cap L|\leq 2
\quad \hbox{ for each } \quad i, \quad 1 \leq i \leq d.
$$

 Therefore, for each $K \subset L'$ or $K \subset L''$
we have $|\ch_i \cap K | \leq 1$, while for $L$ we have
$$
\eta \max_{1 \leq i \leq d} |\ch_i \cap L| \leq {1 \over 2} \max_{1 \leq i \leq d}  |\ch_i \cap L|
\leq \min_{1 \leq i \leq d}  |\ch_i \cap L|.
$$

\medskip

\noindent
{\bf II.2.A.3.} $|\ch \cap L'| < d$ and $|\ch\cap L''| \geq
d$. 
Recall that  $\cc \cap L'\neq \emptyset $ and $\cc\cap L'' \neq \emptyset$, by the defining condition of the case II.2.A.
Then by induction hypothesis all intervals in $\cu \cap L'$  are uncoloured,  but the intervals in $\cc \cap L'$  carry their colours.
Since the pair $(\cc,\cu)$ is $d$-previsible, we have $\cu \cap L'' = \emptyset$.
Therefore, 
$|\ch \cap L''| = |\cc\cap L''|$, and by  condition \eqref{hom2} of the 
$(\eta,d)$-homogeneity for $\cc$, we get
$|\ch_i \cap L''| = |\cc_i \cap L''| \geq 1$ and $|\cc_i \cap L''|$, $1 \leq i \leq d$,  satisfy \rf{(b3)}.

\smallskip

If $\cu \cap L' = \emptyset$ as well, then all intervals in $\ch \cap L$
 come from $\cc \cap L$, and there is nothing to do.

\smallskip

Let $|\cu\cap L'| =x > 0$. We need to colour $x$ intervals in
$\cu\cap L'$. To simplify notation, let $m = |\cc \cap 
L'|$. Note that $1 \leq m+x <d$. Let 
$$S = \{i: |\cc_i \cap L'|
=1\}\quad \hbox{  and  } \quad T= \{i:|\cc_i\cap L'| =0\}.$$
 Let $t_1, \ldots , t_{d-m}$
be an ordering of $T$ such that 
\begin{equation}
\numer{(c4)}
| \cc_{t_1} \cap L''| \leq \ldots \leq |\cc_{t_{d-m}}\cap L''|.
\end{equation}
Since $x < d-m$, there are more colours in $T$ than intervals in $\cU \cap L'$.
Now  attach the colours ${t_1}, \ldots, {t_x}$, bijectively, to intervals in $\cU \cap L'$. 
Then $|\ch_i \cap L'| \leq 1$. 
Consequently, $|\ch_i \cap K| \leq 1$ for each $K \subset L'$.

The colouring of $\ch \cap L'' = \cc \cap L''$ is $(\eta,d)$-homogeneous by the assumption.
It remains  to check that
$|\ch_i \cap L|$ satisfy \rf{(b3)}. Since $\cu \cap L'' = \emptyset$ and $|\ch\cap L''| \geq d$
we have
$$
|\cc \cap L| \geq |\cc \cap L''| = |\ch \cap  L''| \geq d.
$$
Consequently, since the colouring of $\cc$ is $(\eta,d)$-homogeneous, we have $|\cc_i \cap L| \geq |\cc_i \cap L''| \geq 1$ 
and 
\begin{equation}
\numer{(c1)}
\eta \max_{1 \leq i \leq d} |\cc_i \cap  L''|
\leq  \min_{1 \leq i \leq d} |\cc_i \cap  L''|, 
\end{equation}
\begin{equation}
\numer{(c2)}
\eta \max_{1 \leq i \leq d} |\cc_i \cap  L|
\leq  \min_{1 \leq i \leq d}  |\cc_i \cap  L| .
\end{equation}
Moreover,
\begin{equation}
\numer{(c3.a)}
 |\ch_i \cap  L|
=  |\cc_i \cap  L| = |\cc_i \cap  L''|+1 
\text{ for } i \in  S,
\end{equation}
\begin{equation}
\numer{(c3.b)} 
|\ch_i \cap  L|
=  |\cc_i \cap  L| +1 =  |\cc_i \cap  L''|+ 1
 \text{ for } i=t_1, \ldots, t_x,
\end{equation}
\begin{equation}
\numer{(c3.c)}
|\ch_i \cap  L|
= |\cc_i \cap  L| =  |\cc_i \cap  L''|
\text{ for } i=t_{x+1}, \ldots t_{d-m}.
\end{equation}
Let $k$ be such that  $\max_{i} |\ch_i \cap  L| = |\ch_k \cap  L|$.
Then
$k \in S$ or $k \in T$.  

If $k \in S$, then $|\ch_k \cap  L| =|\cc_k \cap  L|$ by \rf{(c3.a)}, and 
\rf{(b3)} is satisfied for $L$ and $\ch$ because of \rf{(c2)} and inequality $|\cc_i\cap L| \leq |\ch_i\cap L|$.

If $k \in T$, then we have 
either $k \in \{t_1, \ldots t_x\}$ or $k \in \{t_{x+1}, \ldots , t_{d-m}\}$.
The ordering defined by \rf{(c4)} implies that  
$k = t_x$ in the first case and 
$k=t_{d-m}$  in the latter case.
If $k=t_{d-m}$, then
\rf{(b3)} is satisfied for $L$ and $\ch$ because of \rf{(c1)}, 
\rf{(c3.c)} and inequality
$|\cc_i\cap L''| \leq |\ch_i \cap L|$.
If $k = t_x$ and $|\ch_{t_x}\cap L| > |\ch_{t_{d-m}}\cap L|$   then we need to check inequality
\begin{equation}
\numer{(c5)}
\eta |\ch_{t_x}\cap L| \leq | \ch_{i}\cap L| \quad {\rm for} \quad 1
\leq i \leq d.
\end{equation}
For $i \in S$  inequality \rf{(c5)} is satisfied because of \rf{(c1)} combined with
\rf{(c3.a)} and \rf{(c3.b)}. For $i= t_1, \ldots, t_x$   inequality \rf{(c5)} is satisfied
because of  \rf{(c2)} and \rf{(c3.b)}. When
$|\ch_{t_x}\cap L| > |\ch_{t_{d-m}}\cap L|$,
 then  \rf{(c3.b)} and \rf{(c3.c)} combined with the ordering \rf{(c4)} 
imply 
$$
|\cc_{t_x} \cap L''| = |\cc_{t_{x+1}} \cap L''|=
 \ldots = |\cc_{t_{d-m}}\cap L''|.
$$ 
 Therefore $|\ch_{t_x} \cap L| = |\cc_{t_x} \cap L''| +1 = |\ch_{i}\cap L| +1$ for all $i=t_{x+1}, \ldots, t_{d-m}$, so
 inequality
 \rf{(c5)} is satisfied, even  with  ${1 \over 2}$ on the left-hand-side, for
 $i=t_{x+1}, \ldots t_{d-m}$.

\bigskip

\noindent
{\bf II.2.A.4.} $|\ch \cap L'| \geq d$ and $|\ch\cap  L''| <
d$. This case is analogous to II.2.A.3.

\medskip

Next, we proceed with the case II.2.B.

\smallskip

\noindent
{\bf II.2.B.1.} $|\ch \cap  L'| \geq d$ and $|\ch\cap  L''| \geq
d$. 
Recall that -- by the condition defining case II.2.B --  $\cc\cap L'=\emptyset$ and $ \cc \cap L''\neq \emptyset$.
Then by the induction hypothesis all intervals in $\cu \cap L''$  are already coloured, but intervals in
$\cu \cap L'$  are uncoloured.

We need to colour intervals in  $\cu \cap L'$. It is enough to
  colour them modulo $d$, see \rf{mod.d}. After this, we get a colouring of $\ch \cap L$ such that both
  $L'$ and $L''$ satisfy \rf{hom2}.
  To check that $L$ satisfies \rf{hom2} as well, we proceed as in case II.2.A.1.

\smallskip

\noindent
{\bf II.2.B.2.} $|\ch \cap  L'| < d$ and $|\ch \cap L''| < d$. 
Induction hypothesis states that  intervals in $\cu \cap L'$ and in $\cu \cap L''$  are uncoloured.
Now we proceed as in case II.2.A.2.

\smallskip

\noindent
{\bf II.2.B.3.} $|\ch \cap L'| < d$ and $|\ch \cap L''| \geq d$.
Recall that $\cc \cap L''\neq \emptyset$. Therefore, by the previsibility assumption, $\cu \cap L'' =\emptyset$.
Since $\cc \cap L' =0$, then by the induction hypothesis intervals in $\cu \cap L'$  are uncoloured, and we need 
to colour them, in case $\cu \cap L' \neq \emptyset$. This is done as in case II.2.A.3.

\smallskip

\noindent
{\bf II.2.B.4.} $|\ch\cap L'| \geq d$ and $|\ch \cap L''| < d$.
Recall that $\cc\cap L'  = \emptyset$. In this case the 
induction hypothesis says that intervals in both $\cu \cap L' $ and $\cu \cap L'' $ 
are uncoloured. We need to colour them all. First, we colour intervals in $\cu \cap L''$  by giving
each of them a different colour which was not used to colour $\cc \cap L''$. This is possible since $|\ch\cap L''| < d$.
Then we colour intervas in $\cu \cap L'$  by modulo $d$ method as in \rf{mod.d}, but starting with colours which have  not been
used to colour intervals in $\ch \cap L''$. In this way we get a modulo $d$ colouring of $\ch \cap L$,  which 
is $({1 \over 2},d)$-homogeneous, hence also $(\eta,d)$-homogeneous.

\smallskip

This completes the proof of Theorem \ref{c.t1}.
\endproof

\section{\numer{c.s3} A colouring problem without solution}

Here we  analyze the role of the  previsibility assumption in Theorem~\ref{c.t1}.

Throughout this section we take 
$d = 2^a$, $a \in \bN$, and $\eta = {1 \over n}$ with $n \in \bN$ and $j \ge n+a +1 .$

We  will define a sequence of collections $\cc(0) \subset \cc(1) \subset \ldots \subset \cc(n) \subset \cd_j$,
of size $|\cc(k)| = k+d$. The initial collection $\cc(0)$ is of size $d$, hence -- up to permutation -- it
 has a unique $(\eta,d)$-homogeneous colouring.
Then we will check that for $1 \leq k \leq n-1$, there is a unique $(\eta,d)$-homogeneous colouring
of $\cc(k)$ keeping the previously determined  $(\eta,d)$-homogeneous colouring
of $\cc(k-1)$. Finally, we will see that there is no $(\eta,d)$-homogeneous colouring
of $\cc(n)$ keeping the previously determined  $(\eta,d)$-homogeneous colouring
of $\cc(n-1)$.

In our example below, the parameter $d$ determines the size of the initial collection $\cc(0)$, while the
parameter $\eta$ determines the number of steps needed to arrive to a problem of consistent colouring without solution.

To define the sequence of collections in question, 
take a chain of dyadic intervals
 $$L_1 \subset L_2 \subset \ldots \subset L_{n+2},
\quad L_i \in \cd_{j-a-i+1}.$$
Then $|L_i| = {1 \over 2} |L_{i+1}|$, and let $P_i$ be the
dyadic brother of $L_i$ in $L_{i+1}$, $ i = 1, \ldots, n+1$. Thus $P_i=L_{i+1} \setminus L_i$.

Now, take two sets of intervals from $\cd_j$:
$$I_1, \ldots, I_{d-1} \in \cd_j \text{ such that } I_i \subset L_1
\text{ for
each } i= 1, \ldots,  d-1,$$
$$J_1, \ldots , J_{n+1} \in \cd_j \text{ such that }
J_i \subset P_i \text{ for each }
i = 1, \ldots , n+1.
$$
Consider the following sequence of collections:

\begin{equation}
\numer{w.k}
\cc(k)  =  \{I_1, \ldots, I_{d-1}\} \cup \{J_{n-k+1}, \ldots, J_{n+1}\}  \hbox{ for } k =0 , \ldots , n.
\end{equation}

\begin{prop}\numer{e.s1} 
The sequence of collection $\cc(k)$, $0 \leq k \leq n$ defined by  \rf{w.k} 
is increasing and it has the following properties:
\begin{itemize}
\item[(A)] Stage $0$.
There exists exactly one -- up to permutation -- $(\eta,d)$-homogeneous  colouring of $\cc(0)$ as
$$
\cc(0) = 
\cc_1(0) \cup \ldots \cup \cc_d(0).
$$

\item[(B)] Stage $k$, $1 \leq k \leq n-1$.
Let 
$$
\cc(k-1) = \cc_1(k-1) \cup \ldots\cup \cc_d(k-1)
$$
be the  $(\eta,d)$-homogeneous colouring of $\cc(k-1)$, obtained at stage $k-1$.
Then there exists exactly one
$(\eta,d)$-homogeneous  colouring of
 $\cc(k)$ as
$$
\cc(k) = 
\cc_1(k) \cup \ldots \cup \cc_d(k),
$$
such that
$$
\cc_i(k-1) \subset \cc_i(k) \hbox{ for all } 1 \leq i \leq d.
$$

\item[(C)] Stage $n$. Let 
$$
\cc(n-1) =  \cc_1(n-1) \cup  \ldots \cup \cc_d(n-1)
$$
be the  $(\eta,d)$-homogeneous colouring of $\cc(n-1)$, obtained at stage $n-1$.
There does not exist an $(\eta,d)$-homogeneous colouring  of  $\cc(n)$ as
$$
\cc(n) = 
\cc_1(n)\cup \ldots\cup \cc_d(n)
$$
such that
$$
\cc_i(n-1) \subset \cc_i(n) \hbox{ for all } 1 \leq i \leq d.
$$
\end{itemize}
\end{prop}

\paragraph{Proof. Verification of (A).}  Consider possible colourings of $\cc(0)$.
Take $L_{n+2}$ as a testing interval. Observe that $|\cc(0)|
= |\cc(0)\cap L_{n+2}|=d$, so  if we want to have $(\eta,d)$-homogeneity, 
we must have \eqref{hom1} and therefore
$ |\cc_i(0)\cap L_{n+2}| =1 $ for each $i=1, \ldots,d$.
 Without loss of generality we can assume
that $J_{n+1}$ has colour 1, and each $I_i$ has colour $i+1$, $i=1, \ldots, d-1$.
Therefore for $\cc(0)$ and each testing interval $L \subset L_{n+2}$ we have
$|\cc_i(0)\cap L|\leq 1$ for each $1 \leq i \leq d$.

\paragraph{The basic observation.} Our example is based on iterating systematically the following basic observation.
Let $k \leq n$.
Assume that $\cc(k)$ has an $(\eta,d)$-homogeneous decomposition as
$$
\cc_1(k), \ldots , \cc_d(k),
$$
so that
$$
\cc_1(0)\subset \cc_1(k), \quad  \ldots , \quad  \cc_d(0) \subset \cc_d(k).
$$
Then necessarily 
\begin{equation}
\numer{jan08.e2}
J_{n-k+1} \text{ must have colour } 1.
\end{equation}
\paragraph{Verification of \eqref{jan08.e2}.}
We know already that  $J_{n+1}$ has to have colour $1$. To check the claim for $J_{n-k+1}$, $k=1, \ldots, n$
we consider the pair of collections $\cc(0) \subset \cc(k)$:
$$
 \cc(k)  =   \cc(0) \cup \{J_{n-k+1}, \ldots,  J_{n}\}.
$$
and testing interval $L_{n-k+2}$. Elements of $\cc(0)$ included in $L_{n-k+2}$ are
$I_1, \ldots, I_{d-1}$. In addition, $J_{n-k+1} \subset P_{n-k+1} \subset L_{n-k+2}$, 
while
$J_{n-k+2}, \ldots, J_n \not \subset L_{n-k+2}$.
Therefore we have
$$
|\cc(0) \cap L_{n-k+2}| = d-1, \quad
|\cc(k) \cap L_{n-k+2}| = d,$$
$$\cc_1(0) \cap L_{n-k+2} = \emptyset \quad
{\rm and} \quad  |\cc_i(0)\cap L_{n-k+2}| = 1 \quad {\rm  for} \quad  i=2, \ldots, d.$$
Therefore, \eqref{hom1} of the $(\eta,d)$-homogeneity condition for 
$\cc(k)$ implies that $J_{n-k+1}$ is of colour 1.

\paragraph{Verificaton of (B).}
Recall that $0 \leq k \leq n$
$$
\cc(k) =  \{I_1, \ldots, I_{d-1}\} \cup \{J_{n-k+1}, \ldots, J_{n+1}\}.
 $$
Moreover, by \eqref{jan08.e2},  the only possible $({1 \over n},d)$-homogeneous decomposition of 
$\cc(k)$ is
$$
\cc_1(k) = \{J_{n-k+1}, \ldots, J_{n+1}\}, \quad \cc_i(k) = \{ I_{i-1} \} \quad \hbox{for} \quad 2 \leq i \leq d
$$
Let's check that for $0 \leq k \leq n-1$, the above decomposition of $\cc(k)$ is indeed $({1 \over n},d)$-homogeneous.
We present the detailed proof for $k=n-1$, since  the cases $ k \leq n-1$ are fully analogous.

First, take as a testing  interval $L_{s}$, $s=3, \ldots n+2$.
Then  elements of $\cc(n-1)$  included in $L_{s}$
are $I_1, \ldots, I_{d-1}$ and $J_2, \ldots, J_{s-1}$.
Therefore $$|\cc(n-1)\cap L_s| = s+d-3,$$ and
$$
|\cc_1(n-1)\cap L_s| = s-2,
\quad |\cc_i(n-1) \cap L_s| = 1 \quad {\rm for} \quad i=2, \ldots d.
$$
Therefore 
$${1 \over n} \max_{1 \leq i \leq d} |\cc_i(n-1)\cap  L_s|
\leq \min_{1 \leq i \leq d} |\cc_i(n-1)\cap L_s|, \quad s = 3, \ldots, n+2.$$

Next take as a testing interval $L_2$.
Then  elements of $\cc(n-1)$  included in $L_{2}$
are $I_1, \ldots, I_{d-1}$, so  $|\cc(n-1) \cap L_2| = d-1$,
 $$
\cc_1(n-1)\cap  L_2 = \emptyset,
\quad |\cc_i(n-1)\cap  L_2| = 1 \quad {\rm for} \quad i=2, \ldots d.
$$
Therefore $L_2$ also satisfies \eqref{hom1}  of the $({1 \over n},d)$-homogeneity condition for $\cc(n-1)$. 
Consequently, $L_1, P_1 \subset L_2$ also satisfy these conditions.

Finally, take as a testing interval $P_k$, $k=2,\ldots, n+1$.
The only element of $\cc(n-1)$ included in $P_k$ is $J_k$,  so $|\cc(n-1)\cap P_k|=1$, and more precisely
 $$
|\cc_1(n-1) \cap P_k| = 1,
\quad \cc_i(n-1) \cap  P_k = \emptyset \quad {\rm for} \quad i=2, \ldots d.
$$
Thus, $P_k$ (and consequently, each testing interval included in $P_k$) satisfies \eqref{hom1} 
of the $({1 \over n},d)$-homogeneity condition for $\cc(n-1)$.

\paragraph{Verification of (C).}
Consider $\cc(n-1)$ and $\cc(n) = \cc(n-1) \cup \{J_1\}$.
Recall that
$$
 \cc(n) =  \{I_1, \ldots, I_{d-1}\} \cup \{J_1, J_2, \ldots, J_{n+1}\}.
$$
Take    $L_{n+2}$ as a testing interval. All intervals from $\cc(n)$ are included in $L_{n+2}$, and the colouring yields
 $$
|\cc_1(n)\cap L_{n+2}| = n+1,
\quad |\cc_i(n)\cap L_{n+2}| = 1 \quad {\rm for} \quad i=2, \ldots d.
$$
For $\cc(n)$ and $L_{n+2}$ we have to consider \eqref{hom2}   of the $({1 \over n},d)$-homogeneity condition.
But the above formulae  mean that for $\cc(n)$ and testing interval $L_{n+2}$, the condition \eqref{hom2}  
 is satisfied
with $\eta' = {1 \over n+1}$, but not with $\eta = {1 \over n}$.\endproof

\paragraph{Remark.} For $0 \leq k \leq n-1$, the pair of collections $(\cc(k), \cu(k))$,
where  $\cU(k) = \cc(k+1) \setminus \cc(k)$ is not $d$-previsible. 
Nevertheless, the colouring problem has a solution for $0 \leq k \leq n-2$.

\bigskip

The examples of Proposition \ref{e.s1} grew out of the couterexamples to classical martingale inequalities, see \cite[p. 105]{stein}, or 
\cite[p. 156]{classical}.

\section{\numer{game} A two-person game}

The problem of consistent colourings
gives rise to the following two-person game. 
The game is played by two players with collections of coloured dyadic intervals in 
$\cD_j$
for a fixed  $j \in \bN$. 
It starts by fixing $\eta > 0,$  $d\in \bN ,$ and a subcollection 
$$
\cC(0) \subset \cD_j
$$
with with an $(\eta,d)$-homogeneous colouring
$$\cc_1(0), \ldots , \cc_d(0),$$
according to Definition \ref{c.d1}. 
The rules of the game are as folows:
\begin{enumerate}
\item In the first stage, Player A chooses a collection $\cC(1) \varsupsetneq \cC(0)$ and $\cC(1) \subset \cD_j$.
Player B determines an $(\eta,d)$-homogeneous colouring of $\cC(1)$ that preserves the colours of $\cC(0)$.
\item At stage $n$, Player A chooses $\cC(n) \varsupsetneq \cC(n-1)$ and $\cC(n) \subset \cD_j$.
Player B determines an $(\eta,d)$-homogeneous colouring of $\cC(n)$  preserving the colours of $\cC(n-1)$.
\item The game stops at stage $n$ if either  $\cC(n-1) = \cD_j$, and then Player B is the winner, 
or else if there does not exist an $(\eta,d)$-homogeneous colouring of $\cC(n)$ that preserves the colours of $\cC(n-1)$.
In the second case, Player A is the winner.
\end{enumerate}

The results of this paper are able to predict  the outcome of the game as follows.
If we do not pose any constraints on the choice of the collections $\cc(k)$, then
the example in Section \ref{c.s3} and Proposition \ref{e.s1} 
describes a realization of our game where Player A has a strategy of winning.

However, if we  restrict the moves of Player A by imposing that 
 $(\cc(k-1), \cc(k) \setminus \cc(k-1))$ is $d$-previsible,  then with the aid of
Theorem \ref{c.t1} and its proof, Player B has always a winning strategy.
In case the moves of Player A  are restricted by $d$-previsibility, we modify the stopping rule
accordingly: Player B is the winner at stage $n$ if there does not exist $\cc(n) \subset \cd_j$ so that $\cc(n) \varsupsetneq \cc(n-1)$
and the pair $(\cc(n-1), \cc(n) \setminus \cc(n-1))$ is $d$-previsible.

\paragraph{Acknowledgement.} A. Kamont was partially supported by SPADE2 programm at IM PAN
and NCN grant N N201 607840. Paul F.X. M\"uller was partially supported by FWF P 20166-N18 and FWF P 23987.

\begin{tabular}{lll}
Anna Kamont & \hspace{4cm} &Paul F. X. M\"uller
\\
Institute of Mathematics & & Department of Analysis
\\
Polish Academy of Sciences & &J.Kepler University
\\
ul. Wita Stwosza 57, 80-952 Gda{\'n}sk & &A-4040 Linz
\\
Poland & & Austria
\\
A.Kamont@impan.gda.pl & & paul.mueller@jku.at
\end{tabular}
\end{document}